\documentclass[12pt,reqno]{amsart}

\usepackage{amssymb,amsmath}

\newcommand{\ov}{\overline}

\newcommand{\N}{{\mathbb N}}
\newcommand{\Z}{{\mathbb Z}}
\newcommand{\Q}{{\mathbb Q}}
\newcommand{\R}{{\mathbb R}}           
\newcommand{\C}{{\mathbb C}}
\newcommand{\HH}{{\mathbb H}}

\newtheorem{theorem}{Theorem}[section]
\newtheorem{lemma}[theorem]{Lemma}
\newtheorem{prop}[theorem]{Proposition}
\newtheorem{cor}[theorem]{Corollary}
\newtheorem*{theoremwithout}{Theorem}

\DeclareMathOperator{\ord}{ord}

\DeclareMathOperator{\tr}{tr}
\DeclareMathOperator{\Tr}{Tr}

\newcommand{\clear}{\setcounter{equation}{0}}

\newcommand{\calD}{{\mathcal D}}

\newcommand{\eps}{\epsilon}

\DeclareMathOperator{\SL}{SL}

\newcommand{\wh}{\widehat}

\allowdisplaybreaks

\begin{document}

\title[]{The Correlation between Multiplicities of \\ closed Geodesics on the Modular Surface}
\author[\relax]{Manfred Peter}
\address{Mathematisches Institut\\
         Albert-Ludwigs-Universit\"at\\
         Eckerstr. 1\\
         D-79104 Freiburg}
\email{peter@arcade.mathematik.uni-freiburg.de}


\date{\today}


\maketitle

\section{Introduction}
\clear

One of the connections between Number Theory and Mathematical Physics that emerged in recent years is 
Arithmetical Quantum Chaos.
On the physical side we have the problem of how the notion of chaos in classical dynamical systems
can be transfered to quantum mechanical systems. Here Gutzwiller's trace formula is a useful quantitative
tool which connects the lengths of closed orbits in the classical picture with the energy eigenvalues
in the quantum mechanical picture. Unfortunately it is an asymptotic relation without rigorous error
estimates. But in the special case of motion on a
surface of constant negative curvature generated by a discontinuous group, Gutzwiller's trace formula is 
Selberg's trace formula and is therefore exact. Furthermore, for particular choices of groups connections
with well understood number theoretical objects can be exploited which are not available for general groups. 

Numerical experiments led to surprising results: For a generic group the eigenvalue statistics of the 
Laplacian on the surface seem to be in accordance with the Gaussian Orthogonal Ensemble (see Bohigas, 
Giannoni, Schmit~\cite{bohigas}). But for arithmetic groups it is closer to Poisson distribution 
(\cite{bohigas}, Aurich, Steiner~\cite{aurich}; see also \cite{bogomolny4}). It seems that the high 
degeneracy of the length spectrum is responsible for the Poissonian behaviour of the eigenvalues. 
Note that a similar phenomenon was observed
by Hejhal and Selberg for quaternion groups (\cite{hejhal2}, Theorems~17.1 and 18.8). They used the high
degeneracy of the length spectrum to prove exceptionally large lower bounds for integral means of the remainder
term in Weyl's asymptotic law for the eigenvalues. Luo and Sarnak~\cite{luo} used the same phenomenon in 
their study of the ``number variance'' (the mean square for the remainder term in short intervals) for
general arithmetic (not necessarily congruence) groups. They proved that for arithmetic groups, the length
spectrum without multiplicities has at most linear growth and conjectured that this property characterizes
arithmetic groups, a fact later proved by Schmutz~\cite{schmutz} in the noncompact case.

In an attempt to understand the Poissonian behaviour for the 
arithmetic group $\Gamma:=\SL(2,\Z)$, Bogomolny, Leyvraz and Schmit~\cite{bogomolny3} calculated the two 
point correlation function for the eigenvalues. Their arguments are not mathematically
rigorous and are given in two steps:
First Selberg's trace formula along with heuristic arguments is used to reduce the pair correlation function
of the eigenvalues to that of the lengths of closed orbits. Second a heuristic version of the 
Hardy-Littlewood method is used to express the latter correlation function as an infinite product with
easily calculable factors. 

Making the first step in \cite{bogomolny3} mathematically
rigorous is extremely hard without a new methodological tool.
Selberg's trace formula seems too weak for this purpose. In the present paper it will be shown how the 
second step can be made rigorous with an approach different from \cite{bogomolny3}. 

There is another way of looking at the result of this paper. Since the structure of the Selberg trace formula
is the same as that of Weil's explicit formula in prime number theory (see~\cite{barner}) one can look upon 
the lengths of closed geodesics as the logarithms of some sort of generalized primes. Thus the theorem
below is an analogue of the Hardy-Littlewood twin prime conjecture. 

In order to state the main result, let $\HH$ be the complex upper half plane and $\Gamma:=\SL(2,\Z)$. This
group acts discontinuously on $\HH$ by M\"obius transformations. Let $F:={}_\Gamma\setminus^\HH$ be the
Riemann surface generated by $\Gamma$. The closed geodesics on $F$ are in one-to-one correspondence to
conjugacy classes of primitive hyperbolic elements in $\Gamma$. For $n\in\N$, $n>2$, let $g(n)$ be the
number of closed geodesics on $F$ which correspond to primitive conjugacy classes with trace $n$. Set
$\alpha(n):=g(n)\log n/n$. Then the function $\alpha$ has mean value $1$. Set 
$\tilde{\alpha}(n):=\alpha(n)-1$.
\begin{theoremwithout} 
  For $r\in \N_0$, the limit
  \[
    \gamma(r):=\lim_{N\to\infty} \frac{1}{N} \sum_{2<n\le N} \tilde{\alpha}(n)\tilde{\alpha}(n+r)
  \]
  exists. Its value is given by 
  \[
    \gamma(r)+1=\prod_p \Big( 1+\sum_{b\ge 1} A_r(p^b) \Big) ,
  \]
  where for a prime $p>2$ and $b\ge 2$, we have
  \begin{align*}
    A_r(p) &= \frac{1}{(p^2-1)^2} \bigg( p\sum_{n\,\bmod\,p} \Big(\frac{(n^2-4)((n+r)^2-4)}{p}\Big)-1
              \bigg),\\
    A_r(p^b) &= \frac{1}{(p^2-1)^2 \, p^{3b-4}} \left\{\begin{array}{ll}
                2p^b(1-1/p) , & r\equiv 0\,(p^b) \\
                -2p^{b-1} , & r\not\equiv 0\,(p^b),\, r\equiv 0\,(p^{b-1}) \\
                \big(\frac{-1}{p}\big)^bp^b(1-1/p) , & r\equiv \pm 4\,(p^b) \\
                -\big(\frac{-1}{p}\big)^bp^{b-1} , & r\not\equiv \pm 4\,(p^b),\, r\equiv\pm 4\, (p^{b-1}) \\
                0 , & \text{otherwise} \end{array}\right\} .
  \end{align*}
  Furthermore, for $b\ge 6$, we have
  \begin{align*}
    & A_r(2) = \frac{1}{9} \left\{\begin{array}{ll} 1, & r\equiv 0\,(2) \\ -1, & r\equiv 1\,(2)
                           \end{array}\right\} ,\quad 
      A_r(4) = \frac{1}{18} \left\{\begin{array}{ll} 1, & r\equiv 0\,(4) \\ -1, & r\equiv 2\,(4) \\
                            0, & r\equiv 1\,(2)  \end{array}\right\} , \\
    & A_r(8)=0,\quad 
      A_r(16) = \frac{1}{9\cdot 16} \left\{\begin{array}{ll} 1, & r\equiv 0\,(16) \\ -1, & r\equiv 8\,(16)\\
                                    0, & r\not\equiv 0\,(8) \end{array}\right\} ,\quad
      A_r(32)=0 , \\
    & A_r(2^b) = \frac{1}{9\cdot 2^{2b-4}} \left\{\begin{array}{ll} 2, & r\equiv 0\,(2^b) \\ 
                                                                -2, & r\equiv 2^{b-1}\,(2^b) \\
                                                                 1, & r\equiv\pm(4+2^{b-2})\,(2^b) \\
                                                                -1, & r\equiv\pm(4+2^{b-2}+2^{b-1})\,(2^b)\\
                                                                 0, & \text{otherwise}
                           \end{array}\right\} .
  \end{align*} 
\end{theoremwithout}
This theorem was already stated in \cite{bogomolny3} and was made plausible on heuristic grounds. In the 
present paper, a different approach is used which exploits the connection of $\alpha(n)$ with a sum of
class numbers of primitive binary quadratic forms (see (\ref{1}) below). The main step is to show that 
$\alpha$ is almost periodic and to compute its Fourier coefficients. Then the theorem follows from 
Parseval's equation. To this end a method is used that has already been applied in \cite{peter2} and 
\cite{peter4}. It might be of interest to note that almost periodic functions --- albeit on the real line 
--- already found applications to trace formulas for integrable geodesic flows (see \cite{bleher7} for
an overview).

The present paper is organized as follows: In section~2 the function $\alpha$ is reduced to class numbers.
Section~3 contains a short review of almost periodic arithmetical functions. Furthermore, the function
$\alpha$ is shown to be almost periodic by approximating it with a carefully choosen periodic function.
This will considerably simplify the computation of the product representation of $\gamma(r)$. In section~4
this Euler product is derived and its local factors are computed in section~5. 

{\bf Acknowledgement.} I would like to express my sincere gratitude to Prof. Ze\'ev Rudnick for bringing this
problem to my attention. I gained much from conversations with him and from the stimulating atmosphere
which he and his co-organizers created at the DMV seminar ``The Riemann Zeta Function and Random Matrix
Theory''.

\section{Reduction to class numbers}
\clear

In order to state the connection with class numbers we need some notation. Proofs of the number theoretic
facts in this section can be found for example in \cite{borevich}, \cite{davenport} and \cite{ono}.

The letter $d$ will always stand for a positive non-square discriminant (i.e. $d\equiv 0,1\,(4)$). A
primitive binary quadratic form is a polynomial $f=ax^2+bxy+cy^2$ with $a,b,c\in\Z$ and $\gcd(a,b,c)=1$.
Its discriminant is $d=b^2-4ac$. Two such forms $f_i$, $i=1,2$, are called equivalent if there is a matrix
$\begin{pmatrix} \alpha & \beta \\ \gamma & \delta \end{pmatrix} \in \SL(2,\Z)$ with 
$f_1(x,y)=f_2(\alpha x+\beta y,\gamma x+\delta y)$. A fundamental theorem in number theory now states that
the number $h(d)$ of equivalence classes of primitive binary quadratic forms with discriminant $d$ is finite. This
number is one of the important quantities in number theory since it appears in a surprising variety of 
situations.

The Pellian equation $u^2-dv^2=4$ has infinitely many solutions in integers $(u,v)$. The solution 
$(u_d,v_d)$ with $u_d,v_d>0$ and $u_d$ minimal is called the fundamental solution since all the other 
solutions can be generated from it by a simple law. Set $\eps_d:=(u_d+v_d\sqrt{d})/2$. The trace of $\eps_d$
is defined as $\tr(\eps_d):=u_d$. 

The one-to-one correspondence between primitive conjugacy classes in $\Gamma$ and closed geodesics on $F$
is as follows: Every primitive hyperbolic $P\in\Gamma$ has two real fixed points (one of them can be 
$\infty$). The orthogonal circle in $\HH$ which ends in these fixed points induces a closed geodesic
on $F$ of length $l$ where $\cosh(l/2)=|\Tr P|$. The possible values for $l$ and their multiplicities
are described in 
\begin{prop}[\cite{sarnak}, Corollary~1.5]
  The lengths of closed geodesics on $F$ are the numbers $2\log\eps_d$ with multiplicities $h(d)$.
\end{prop}
From this proposition and the preceding description it follows that 
\begin{equation} \label{1}
  \alpha(n) = \frac{\log n}{n} \sum_{d:\,\tr(\eps_d)=n} h(d) .
\end{equation}
Quantitative results about class numbers are often derived using Dirichlet's class number formula. First we
must define Jacobi's character and Dirichlet L-series. Define $\chi_d:\Z\to\{0,\pm1\}$ to be completely
multiplicative with 
\begin{align*}
  & \chi_d(p):= \left\{\begin{array}{ll} 1, & x^2\equiv d\,(p)\,\text{is solvable},\, p\not| d \\
                                        -1, & x^2\equiv d\,(p)\,\text{is insolvable} \\
                                         0, & p|d \end{array}\right\} , \\ 
  & \chi_d(2):= \left\{\begin{array}{ll} 1, & d\equiv 1\,(8)\\
                                        -1, & d\equiv 5\,(8) \\
                                         0, & d\equiv 0\,(4) \end{array}\right\} , \quad
    \chi_d(-1):=1 .
\end{align*}
$\chi_d$ is a Dirichlet character modulo $d$ (the quadratic reciprocity law is used to prove the
$d$-periodicity). For $\Re s>0$, the series 
\[
  L(s,\chi_d) := \sum_{n\ge 1} \frac{\chi_d(n)}{n^s}
\]
is uniformly convergent and defines a holomorphic function. 
\begin{prop}[Dirichlet's class number formula] \label{101}
  For a positive non-square discriminant $d$, we have
  \[
    h(d)\log\eps_d = \sqrt{d}\, L(1,\chi_d) .
  \]
\end{prop}
In the next section this formula will be used to prove the almost periodicity of $\alpha$. In principle 
this could be done by writing 
\[
  L(1,\chi_d) = \sum_{1\le n\le N} \frac{\chi_d(n)}{n} + \text{error} .
\]
Note that the sum on the right hand side is a periodic function of $d$ and can be used to approximate
$L(1,\chi_d)$. But this procedure has three severe drawbacks: First the approximation is not particularly
good. Therefore we will use a smoothed version of the series representation of $L(1,\chi_d)$ instead.
Second we must approximate a sum of values $L(1,\chi_d)$ with the condition 
$\tr(\eps_d)=n$. Breaking up this condition into easier summations, as must be done to process it further,
would make things much more difficult. Third we want to compute the Fourier coefficients of 
$\alpha$ and show that they are multiplicative. This would be near to impossible with the above 
approach. 

Instead an approximating periodic function is used which already incorporates some sort of 
multiplicativity. From the multiplicativity of $\chi_d$ the Euler product 
\[
  L(s,\chi_d) = \prod_p \Big( 1-\frac{\chi_d(p)}{p^s} \Big)^{-1} ,\quad \Re s>1 ,
\]
follows. Thus it seems more reasonable to use a partial product of this representation than a partial sum 
of the series as approximating function. 

\section{Almost periodicity}
\clear

The standard reference for almost periodic arithmetical functions is \cite{schwarz}. Here the necessary
material will be reviewed briefly. 

Let $q\ge 1$. For $f:\N\to\C$, define the seminorm 
\[
  \|f\|_q := \bigg( \limsup_{N\to\infty} \frac{1}{N}\sum_{1\le n\le N} |f(n)|^q \bigg)^{1/q} \in
  [0,\infty] .
\]
$f$ is called $q$-limit periodic if for every $\eps>0$ there is a periodic function $h$ with 
$\|f-h\|_q\le\eps$. The set $\calD^q$ of all $q$-limit periodic functions becomes a Banach space with
norm $\|\cdot\|_q$ if functions $f_1,f_2$ with $\|f_1-f_2\|_q=0$ are identified. If 
$1\le q_1\le q_2<\infty$, we have $\calD^1\supseteq \calD^{q_1}\supseteq\calD^{q_2}$ as sets (but they 
are endowed with different norms!). There is the more general notion of $q$-almost periodic function
which will not be used in this paper (they are defined as above but with arbitrary trigonometric polynomials 
for $h$ instead of periodic functions). For all $f\in\calD^1$, the mean value 
\[
  M(f):=\lim_{N\to\infty} \frac{1}{N}\sum_{1\le n\le N} f(n)
\]
exists. The space $\calD^2$ is a Hilbert space with inner product 
\[
  \langle f,h\rangle := M(f\ov{h}) ,\quad f,h\in\calD^2 .
\]
For $u\in\R$, define $e_u(n):=e^{2\pi iun}$, $n\in\N$. For all $f\in\calD^1$, the Fourier coefficients 
$\wh{f}(u):=M(fe_{-u})$, $u\in\R$, exist. For $u\not\in\Q$, we have $\wh{f}(u)=0$ (this comes from the
fact that $f$ can be approximated by linear combinations of functions $e_v$ with $v\in\Q$, and that
$M(e_ve_{-u})$ equals $1$ if $v-u\in\Z$ and $0$ otherwise; for almost periodic functions it is no longer 
true). In $\calD^2$, we have the canonical orthonormal base $\{e_{a/b}\}$, where $1\le a\le b$ and 
$\gcd(a,b)=1$. 

Limit periodic (and almost periodic) functions have a couple of nice properties. They can be added,
multiplied and plugged into continuous functions and, under certain conditions, the result is again a
limit periodic (almost periodic) function. They have mean values and limit distributions. Here we will 
use Parseval's
equation. In Corollary~\ref{109} we will prove that $\alpha\in\calD^q$ for all $q\ge 1$. As a side result
in section~5, we get $M(\alpha)=\wh{\alpha}(0)=1$. Thus for $r\in\N_0$, we have $\tilde{\alpha},
\tilde{\alpha}_{+r}\in\calD^2$ (where $\tilde{\alpha}_{+r}(n):=\tilde{\alpha}(n+r)$), and Parseval's 
equation gives 
\begin{align}
  \gamma(r) &= M(\tilde{\alpha}\ov{\tilde{\alpha}_{+r}}) = \langle \alpha,\alpha_{+r}\rangle
  -2M(\alpha)+1 \nonumber \\
  & = \sum_{b\ge 1}\sum_{1\le a\le b:\, (a,b)=1} \wh{\alpha}\Big(\frac{a}{b}\Big)
  \ov{\wh{\alpha_{+r}}\Big(\frac{a}{b}\Big)} -1 
  = \sum_{b\ge 1} \sum_{1\le a\le b:\, (a,b)=1} \Big|\wh{\alpha}\Big(\frac{a}{b}\Big)\Big|^2e^{2\pi i
  ar/b} -1 . \label{2}
\end{align}
The last equation follows easily since $\alpha$ is real valued. If we can compute $\wh{\alpha}$ and show
that it is multiplicative our main theorem follows. 

The approximation of $\alpha$ by periodic functions is done in two steps. First we bring Dirichlet 
L-series into the picture. For $n\in\N$, $n>2$, define
\[
  \beta(n):=\sum_{d,v\ge 1:\, dv^2=n^2-4} \frac{1}{v}L(1,\chi_d) ,
\]
where we must remember that $d$ will always run through non-square positive discriminants. 
\begin{lemma} \label{102}
  For $q\ge 1$, we have $\|\alpha-\beta\|_q=0$.
\end{lemma}
\begin{proof}
For $d$ fixed, the powers $\eps_d^l$, $l\ge 1$, of the fundamental unit give all solutions $(u,v)$ of the
Pellian equation $u^2-dv^2=4$ with $u,v\ge 1$ by way of the rule
\[
  \frac{u+v\sqrt{d}}{2} = \eps_d^l .
\]
For every such solution we have $v\sqrt{d}=\sqrt{u^2-4}$. Thus Proposition~\ref{101} gives 
\[
  \beta(n)=\sum_{d,l\ge 1:\, \tr(\eps_d^l)=n} \frac{h(d)\log\eps_d}{\sqrt{n^2-4}} .
\]
Since 
\[
  \log\eps_d=\frac{1}{l}\log\eps_d^l=\frac{1}{l}\log\frac{n+\sqrt{n^2-4}}{2} ,
\]
it follows from (\ref{1}) that 
\begin{align*}
  \beta(n) &= \frac{1}{\sqrt{n^2-4}}\log\frac{n+\sqrt{n^2-4}}{2} \sum_{d,l\ge 1:\,\tr(\eps_d^l)=n}
  \frac{1}{l}h(d) \\
  &= \frac{n}{\sqrt{n^2-4}} \frac{\log\frac{1}{2}(n+\sqrt{n^2-4})}{\log n} \alpha(n)
  +O\bigg( \frac{\log n}{n}\sum_{d\ge 1,\, l\ge 2:\, \tr(\eps_d^l)=n} h(d) \bigg) .
\end{align*}
Since $h(d)\ll\sqrt{d}$ (this can be seen, e.g. from Proposition~\ref{101} and the estimates 
$\log\eps_d\ge \log\frac{1}{2}(1+\sqrt{d})$ and $L(1,\chi_d)\ll\log d$; the latter follows easily by 
partial summation from the orthogonality relation for characters), it follows that for $d,l\ge 1$
with $\tr(\eps_d^l)=n$,
\[
  h(d)\ll\sqrt{d}\ll(\eps_d^l)^{1/l}=\Big(\frac{n+\sqrt{n^2-4}}{2}\Big)^{1/l}\le n^{1/l} .
\]
Thus 
\[
  \alpha(n) \ll \log n\sum_{d,v\ge 1:\, dv^2=n^2-4} 1 \le \log n \,\,\tau(n^2-4) \ll_\eps n^\eps
\]
for all $\eps>0$; here $\tau$ denotes the divisor function and $\tau(m)\ll_\eps m^\eps$ was used. This
implies
\begin{align*}
  |\beta(n)&-\alpha(n)| \\
  &\ll \bigg|\Big(1-\frac{4}{n^2}\Big)^{-1/2}\bigg(1+\frac{\log\frac{1}{2}
  (1+\sqrt{1-4/n^2})}{\log n}\bigg)-1\bigg| \cdot |\alpha(n)| + \frac{\log n}{\sqrt{n}} \sum_{d,v\ge 1:\,
  dv^2=n^2-4} 1 \\
  & \ll n^{-1/2+\eps} ,
\end{align*}
which proves the lemma.
\end{proof}
In the crucial second step $\beta$ is approximated by a sum of partial products of the Euler product of
$L(1,\chi_d)$. For $n>2$, $P\ge 2$, define
\[
  \beta_P(n):=\sum_{d,v\ge 1:\, dv^2=n^2-4,\, p|v \Rightarrow p\le P} \frac{1}{v}
  \prod_{p\le P} \Big(1-\frac{\chi_d(p)}{p}\Big)^{-1} .
\]
Then $\beta(n)-\beta_P(n)=\Delta_P^{(1)}(n)+\Delta_P^{(2)}(n)$, where
\[
  \Delta_P^{(1)}(n) := \sum_{d,v\ge 1:\, dv^2=n^2-4 ,\, p|v \text{ for some } p>P} \frac{1}{v}
  L(1,\chi_d)
\]
and
\[
  \Delta_P^{(2)}(n) := \sum_{d,v\ge 1:\, dv^2=n^2-4 ,\, p|v \Rightarrow p\le P} \frac{1}{v}
  \bigg( L(1,\chi_d)-\prod_{p\le P}\Big(1-\frac{\chi_d(p)}{p}\Big)^{-1} \bigg) .
\]
Fix $q\in\N$ with $q>1$ and choose $q'>1$ with $1/(2q)+1/q'=1$. The next lemma shows that
$\Delta_P^{(1)}(n)$ is negligible in the $2q$-norm as $P\to\infty$. 
\begin{lemma} \label{103}
  For $P\ge 2$, we have
  \[
    \|\Delta_P^{(1)}\|_{2q} \ll \bigg(\sum_{v>P}\frac{1}{v^{q'}} \bigg)^{1/q'} .
  \]
\end{lemma}
\begin{proof}
H\"older's inequality gives 
\[
  |\Delta_P^{(1)}(n)| \le \bigg( \sum_{d,v\ge 1:\, dv^2=n^2-4 \atop p|v \text{ for some } p>P}
  \frac{1}{v^{q'}} \bigg)^{1/q'} \bigg( \sum_{d,v\ge 1:\, dv^2=n^2-4 \atop p|v \text{ for some } p>P}
  L(1,\chi_d)^{2q} \bigg)^{1/(2q)} .
\]
For $x\ge 1$, this gives
\[
  \sum_{2<n\le x} |\Delta_P^{(1)}(n)|^{2q} \le \bigg( \sum_{v>P} \frac{1}{v^{q'}} \bigg)^{2q/q'}
  \sum_{2<n\le x,\, d,v\ge 1:\, dv^2=n^2-4} L(1,\chi_d)^{2q} .
\]
The second sum on the right hand side is 
\[
  \sum_{d,l\ge 1:\, \tr(\eps_d^l)\le x} L(1,\chi_d)^{2q} \sim \text{const}\cdot x
\]
as $x\to\infty$ (see \cite{peter3}). This means that the values $L(1,\chi_d)$ are constant in the mean
when ordered according to the sizes of their fundamental units. Thus the lemma follows. 
\end{proof}
In order to estimate $\Delta_P^{(2)}(n)$ we must compare $L(1,\chi_d)$ with a partial product of its Euler
product. This is done by comparing both terms with a smoothed version of the Dirichlet series for
$L(1,\chi_d)$. Let $N\ge 1$. Then 
\[
  \Delta_P^{(2)}(n) = \Delta_{P,N}^{(2,1)}(n) + \Delta_{P,N}^{(2,2)}(n) + \Delta_{P,N}^{(2,3)}(n) ,
\]
where
\begin{align*}
  \Delta_{P,N}^{(2,1)}(n) &:= \sum_{d,v} \frac{1}{v} \bigg( L(1,\chi_d) - \sum_{l\ge 1} 
  \frac{\chi_d(l)}{l}\,\,e^{-l/N} \bigg) , \\
  \Delta_{P,N}^{(2,2)}(n) &:= \sum_{d,v} \frac{1}{v} \sum_{l\ge 1:\, p|l \text{ for some } p>P}
  \frac{\chi_d(l)}{l} \,\,e^{-l/N} , \\
  \Delta_{P,N}^{(2,3)}(n) &:= \sum_{d,v} \frac{1}{v} \sum_{l\ge 1:\, p|l \Rightarrow p\le P}
  \frac{\chi_d(l)}{l} \Big(e^{-l/N}-1\Big) ;
\end{align*}
here the conditions on $d$ and $v$ are as in $\Delta_{P}^{(2)}(n)$. The third term can be estimated
easily.
\begin{lemma} \label{104}
  For $P\ge2$ and $x,N\ge 1$, we have
  \[
    \bigg( \frac{1}{x} \sum_{2<n\le x} \big|\Delta_{P,N}^{(2,3)}(n)\big|^{2q} \bigg)^{1/(2q)}
    \ll N^{-1/2} + \sum_{l>\sqrt{N}:\, p|l \Rightarrow p\le P} \frac{1}{l} .
  \]
\end{lemma}
\begin{proof}
Since $|e^{-u}-1|\ll u$ for $0\le u\le 1$, we see that for $n>2$ the inner sum in 
$\Delta_{P,N}^{(2,3)}(n)$ is
\begin{align*}
  & \ll \sum_{l\ge 1:\, p|l \Rightarrow p\le P} \frac{1}{l} \big| e^{-l/N}-1\big| 
  \ll \sum_{l>\sqrt{N}:\, p|l \Rightarrow p\le P} \frac{2}{l} + \sum_{1\le l\le \sqrt{N}} \frac{1}{l}\cdot
    \frac{l}{N} \\
  & \ll \sum_{l>\sqrt{N}:\, p|l \Rightarrow p\le P} \frac{1}{l} + N^{-1/2} =: c_1(P,N) .
\end{align*}
H\"older's inequality now gives
\begin{align*}
  \sum_{2<n\le x} \big|\Delta_{P,N}^{(2,3)}(n)\big|^{2q} &\ll \sum_{2<n\le x}
    \bigg(\sum_{d,v\ge 1:\, dv^2=n^2-4} \frac{1}{v} \bigg)^{2q} c_1(P,N)^{2q} \\
  & \ll c_1(P,N)^{2q} \sum_{2<n\le x} \bigg( \sum_{d,v\ge 1:\, dv^2=n^2-4} \frac{1}{v^{q'}} \bigg)^{2q/q'}
    \bigg( \sum_{d,v\ge 1:\, dv^2=n^2-4} 1 \bigg) \\
  & \ll c_1(P,N)^{2q} \sum_{2<n\le x,\, d,v\ge 1:\, n^2-dv^2=4} 1 ,
\end{align*}
where $q'>1$ was used. The last sum is 
\[
  \sum_{d,l\ge 1:\, \tr(\eps_d^l)\le x} 1 \sim \text{const} \cdot x 
\]
(see \cite{peter3} or \cite{sarnak}), and the lemma follows. 
\end{proof}
So far we did not use the oscillation of the Jacobi character. For the estimation of 
$\Delta_{P,N}^{(2,2)}$ we must take it into account.
\begin{lemma} \label{105}
  For $l,v\in\N$ and $x\ge 3$, we have
  \[
    \sum_{d\ge 1,\, 2<n\le x:\, dv^2=n^2-4} \chi_d(l) \ll_\eps \frac{x}{v^{2-\eps}K(l)} + v^\eps l ,
  \]
  where $K(l)$ is the squarefree kernel of $l$ and $\eps>0$.
\end{lemma}
\begin{proof} See \cite{peter4}, estimate~(2.7). \end{proof}
\begin{lemma} \label{106}
  For $P\ge 2$ and $x,N\ge 1$, we have
  \[
    \frac{1}{x}\sum_{2<n\le x} \big|\Delta_{P,N}^{(2,2)}(n)\big|^{2q} \ll \sum_{l>P^{2q}} 
    \frac{\tau_{2q}(l)}{l\,K(l)} + \frac{1}{x} (\log N)^{2q} + x^{-1/2+\eps} N^{2q} ,
  \]
  where $\tau_k(l):=\sum_{l_1,\ldots, l_k\ge 1:\, l_1\cdots l_k=l} 1$.
\end{lemma}
\begin{proof}
Split the first sum in $\Delta_{P,N}^{(2,2)}(n)$ into two sums depending on whether $v\le \sqrt{n}$
or $v>\sqrt{n}$. Thus $\Delta_{P,N}^{(2,2)}(n)=\Delta_{P,N}^{(2,2,1)}(n)+\Delta_{P,N}^{(2,2,2)}(n)$. 
A trivial estimate gives
\[
  \Delta_{P,N}^{(2,2,2)}(n) \ll \sum_{d,v\ge 1:\, dv^2=n^2-4,\, v>\sqrt{n}}\,\, \frac{1}{v} \sum_{l\ge 1}
  \frac{1}{l} \,\,e^{-l/N} \ll \log N \frac{1}{\sqrt{n}}\,\,\tau(n^2-4)
\]
and thus 
\begin{equation} \label{3}
  \sum_{2<n\le x} \big|\Delta_{P,N}^{(2,2,2)}(n)\big|^{2q} \ll_\eps (\log N)^{2q} \sum_{2<n\le x}
  n^{2q(-1/2+\eps)} \ll_\eps (\log N)^{2q}
\end{equation}
since $q>1$. H\"older's inequality gives 
\begin{align*}
  & \big|\Delta_{P,N}^{(2,2,1)}(n)\big| \le \sum_{d,v\ge 1:\, dv^2=n^2-4,\, v\le \sqrt{n}} \,\, \frac{1}{v}
    \bigg|\sum_{l\ge 1:\, p|l \text{ for some } p>P} \frac{\chi_d(l)}{l} \,\,e^{-l/N} \bigg| \\
  & \le \bigg( \sum_{d,v\ge 1:\, dv^2=n^2-4} \frac{1}{v^{q'}} \bigg)^{1/q'} \bigg( 
    \sum_{d,v\ge 1:\, dv^2=n^2-4,\, v\le \sqrt{n}} \Big( \sum_{l\ge 1:\, p|l \text{ for some } p>P}
    \frac{\chi_d(l)}{l}\,\,e^{-l/N} \Big)^{2q} \bigg)^{1/(2q)} .
\end{align*}
Thus for $x\ge 1$, 
\begin{align*}
  & \sum_{2<n\le x} \big|\Delta_{P,N}^{(2,2,1)}(n)\big|^{2q} \ll \sum_{2<n\le x}\,\,\, \sum_{d,v\ge 1:\, dv^2=
    n^2-4,\, v\le\sqrt{x}} \bigg( \sum_{l\ge 1:\, p|l \text{ for some } p>P} \frac{\chi_d(l)}{l}
    \,\,e^{-l/N} \bigg)^{2q} \\
  & = \sum_{l_1,\ldots, l_{2q}:\, p_i|l_i \text{ for some } p_i>P} \frac{1}{l_1\cdots l_{2q}}\,\,
    e^{-(l_1+\cdots +l_{2q})/N} \sum_{1\le v\le \sqrt{x}}\,\,\, \sum_{2<n\le x,\, d\ge 1:\, dv^2=n^2-4}
    \chi_d(l_1\cdots l_{2q}) .
\end{align*}
Applying Lemma~\ref{105} to the innermost sum gives the estimate
\begin{align*}
  & \ll \sum_{l>P^{2q}} \frac{1}{l}\,\, \tau_{2q}(l) \sum_{1\le v\le \sqrt{x}} \frac{x}{v^{2-\eps}K(l)}
    + \sum_{l_1,\ldots, l_{2q}\ge 1} \frac{l_1\cdots l_{2q}}{l_1\cdots l_{2q}}\,\, e^{-(l_1+\cdots+l_{2q})/N}
    \sum_{1\le v\le \sqrt{x}} v^\eps \\
  & \ll x\sum_{l>P^{2q}} \frac{\tau_{2q}(l)}{l\,K(l)} + x^{(1+\eps)/2} N^{2q} ,
\end{align*}
which together with (\ref{3}) proves the lemma.
\end{proof}
In order to estimate $\Delta_{P,N}^{(2,1)}(n)$ we must show that the error 
\begin{equation} \label{10}
  I(d,N) := L(1,\chi_d) - \sum_{l\ge 1} \frac{\chi_d(l)}{l} \,\,e^{-l/N} ,
\end{equation}
which comes from smoothing the Dirichlet series expansion of $L(1,\chi_d)$, is small for large $N$.
This is done by representing $I(d,N)$ as an integral over a vertical line in the critical strip 
$\{0<\Re s<1\}$ and using information about the location of the non-trivial zeros of $L(s,\chi_d)$. 
The Dirichlet series for $L(s,\chi_d)$ is absolutely and uniformly convergent on the line $\Re s=2$.
Sterling's formula gives
\[
  \Gamma(s) \ll |\Im s|^c \,\,e^{-\pi|\Im s|/2}
\]
for $c_1\le \Re s\le c_2$, $|\Im s|\ge 1$, with some constant $c=c(c_1,c_2)$. Using Mellin's formula
\[
  \frac{1}{2\pi i} \int_{1-i\infty}^{1+i\infty} \Gamma(s)y^{-s} ds = e^{-y} ,\quad y>0 ,
\]
we get 
\begin{equation} \label{4}
  \frac{1}{2\pi i} \int_{2-i\infty}^{2+i\infty} \Gamma(s-1)\,L(s,\chi_d)\,N^{s-1}ds = \sum_{l\ge 1}
  \frac{\chi_d(l)}{l} \,\,e^{-l/N} .
\end{equation}
On the other hand, for $1/2\le \Re s\le 2$, we have
\begin{equation} \label{5}
  L(s,\chi_d) \ll (d|s|)^{1/2}
\end{equation}
(this is easily seen by partial summation). Thus the line of integration in (\ref{4}) may be moved to the
line $\Re s=\kappa$ with some $1/2<\kappa<1$. Taking into account the pole of the integrand at $s=1$,
and using the residue theorem we get for (\ref{4}) the expression
\[
  \frac{1}{2\pi i} \int_{\kappa-i\infty}^{\kappa+i\infty} \Gamma(s-1)\,L(s,\chi_d)\,N^{s-1}ds + L(1,\chi_d)
\]
and thus 
\begin{equation} \label{6}
  I(d,N) = -\frac{1}{2\pi i} \int_{\kappa-i\infty}^{\kappa+i\infty} \Gamma(s-1)\,L(s,\chi_d)\,N^{s-1}ds .
\end{equation}
We must now considerably reduce the exponent $1/2$ of $d$ in (\ref{5}). In order to see the principle let us
first assume the Generalized Riemann Hypothesis which says that for all Dirichlet characters $\chi$ modulo
$q$, the L-series $L(s,\chi)$ has only zeros with real part $1/2$ in the critical strip $0<\Re s<1$. From
this the Generalized Lindel\"of Hypothesis follows and in particular for all $d\ge 1$, $\Re s=\kappa$ 
and $\eps>0$, we have 
\[
  L(s,\chi_d) \ll_\eps (d|s|)^\eps .
\]
From (\ref{6}) it follows that for $d,N\ge 1$, 
\[
  I(d,N) \ll_\eps d^\eps N^{\kappa-1} .
\]
This shows that 
\[
  \Delta_{P,N}^{(2,1)}(n) \ll \sum_{d,v\ge 1:\, dv^2=n^2-4} |I(d,N)| \ll_\eps N^{\kappa-1} n^{3\eps}
\]
and for $P\ge 2$ and $x,N\ge 1$, we get 
\begin{equation} \label{7}
  \frac{1}{x} \sum_{2<n\le x} \big|\Delta_{P,N}^{(2,1)}(n)\big|^{2q} \ll_\eps x^{6q\eps} N^{2q(\kappa-1)} .
\end{equation}
Here it is important that the exponent of $N$ is negative. Taking $N$ to be a power of $x$ with small 
exponent therefore lets the contribution of $\Delta_{P,N}^{(2,1)}(n)$ vanish as $x\to\infty$. The next lemma
gives an estimate which for our purposes is as good as (\ref{7}) and can be proved unconditionally. 
\begin{lemma} \label{107}
  There are $0<\kappa,\mu<1$ such that for $P\ge 2$, $x,N\ge 1$ and $\eps>0$, we have
  \[
    \frac{1}{x}\sum_{2<n\le x} \big|\Delta_{P,N}^{(2,1)}(n)\big|^{2q} \ll_\eps x^\eps N^{2q(\kappa-1)}
    + x^{\mu-1+\eps} \big(\log(x^2N)\big)^{2q} .
  \]
\end{lemma}
\begin{proof}
Choose $1/2<\sigma_0<1$ with $\mu:=8(1-\sigma_0)/\sigma_0<1$. Choose $\sigma_0<\kappa<1$. Define the 
rectangle
\[
  R_x := \big\{ s\in\C \,\big|\, \sigma_0\le \Re s\le 1,\, |\Im s|\le \log^2 x \big\} .
\]
If  $d\le x^2$ and $L(s,\chi_d)$ has no zeros in $R_x$ then a standard argument (see, for example, 
Titchmarsh~\cite{titchmarsh}, Theorem~14.2) shows that for $\Re s=\kappa$, $|\Im s|\le (\log x)^2/2$, we
have
\[
  \log|L(s,\chi_d)| \ll \log\log(x+2)\,\big(\log(x+2)\big)^{(1-\kappa)/(1-\sigma_0)} \le \eps\log x
  +c(\eps) 
\]
with some constant $c(\eps)>0$ depending on $\eps$. Together with (\ref{5}) and (\ref{6}) this gives 
\begin{align}
  I(d,N) &\ll_\eps \int_{|t|\le (\log x)^2/2} (|t|+1)^c\,\,e^{-\pi|t|/2}\,x^\eps\,N^{\kappa-1}dt \nonumber \\
  & \hspace{0.7cm} + \int_{|t|\ge (\log x)^2/2} |t|^c\,\,e^{-\pi|t|/2}\, (d|t|)^{1/2}\,N^{\kappa-1} dt 
    \nonumber \\
  & \ll_\eps x^\eps\,N^{\kappa-1} .  \label{8}
\end{align} 
Next we must show that $L(s,\chi_d)$ cannot have a zero in $R_x$ too often. From zero density estimates 
it follows that 
\begin{equation} \label{9}
  \#\big\{ (n,v,d) \,\big|\, 2<n\le x,\, d,v\ge 1,\, n^2-dv^2=4,\, L(s,\chi_d) \text{ has a zero in } R_x
  \big\} \ll x^{\mu+\eps}
\end{equation}
(see \cite{peter3}, Lemma~4.11 or \cite{peter4}, estimate~(2.6)). A trivial estimation of (\ref{10})
gives 
\begin{equation} \label{11}
  I(d,N) \ll \log (dN) .
\end{equation}
(\ref{8}), (\ref{9}), (\ref{11}) and H\"older's inequality give 
\begin{align*}
  & \sum_{2<n\le x} \big|\Delta_{P,N}^{(2,1)}(n)\big|^{2q} \ll \sum_{2<n\le x} \bigg( \sum_{d,v\ge 1:\,
    dv^2=n^2-4} \frac{1}{v^{q'}} \bigg)^{2q/q'} \bigg( \sum_{d,v\ge 1:\, dv^2=n^2-4} |I(d,N)|^{2q} \bigg) \\
  & \ll \sum_{2<n\le x,\, d,v\ge 1:\, dv^2=n^2-4 \atop L(s,\chi_d) \text{ has no zeros in } R_x} 
    (x^\eps\,N^{\kappa-1})^{2q} 
    + \sum_{2<n\le x,\, d,v\ge 1:\, dv^2=n^2-4 \atop L(s,\chi_d) \text{ has a zero in } R_x}
    \big(\log (dN)\big)^{2q} \\
  & \ll x(x^\eps N^{\kappa-1})^{2q} + x^{\mu+\eps} \big(\log(x^2N)\big)^{2q} ,
\end{align*}
which proves the lemma.
\end{proof}
Now the results are collected.
\begin{prop} \label{108}
  For $P\ge 2$, we have 
  \[
    \|\beta-\beta_P\|_{2q} \ll \bigg(\sum_{v>P} \frac{1}{v^{q'}}\bigg)^{1/q'} + \bigg( \sum_{l>P^{2q}}
    \frac{\tau_{2q}(l)}{l\,K(l)} \bigg)^{1/(2q)} .
  \]
\end{prop}
\begin{proof}
For $x\ge 1$, choose $N:=x^{1/(8q)}$. Then Lemmas~\ref{104}, \ref{106} and \ref{107} show that
\begin{align*}
  \frac{1}{x} \sum_{2<n\le x} \big|\Delta_P^{(2)}(n)\big|^{2q} \ll & \bigg( x^{-1/(16q)} + 
    \sum_{l>x^{1/(16q)}:\, p|l \Rightarrow p\le P} \frac{1}{l} \bigg)^{2q} 
  + \sum_{l>P^{2q}} \frac{\tau_{2q}(l)}{l\,K(l)} \\
  & + \frac{1}{x}(\log x)^{2q} + x^{-1/4+\eps}
    +x^{(\kappa-1)/4+\eps} + x^{\mu-1+\eps}(\log x)^{2q} .
\end{align*}
Since the series $\sum_{l\ge 1:\, p|l \Rightarrow p\le P} 1/l$ converges, we have for $P\ge 2$ fixed
\[
  \|\Delta_P^{(2)}(n)\|_{2q}^{2q} \ll \sum_{l>P^{2q}} \frac{\tau_{2q}(l)}{l\,K(l)} .
\]
Together with Lemma~\ref{103} this proves the proposition.
\end{proof}

\section{An euler product}
\clear

In this section we exploit the particular construction of $\beta_P$ by writing it as a product of functions
each depending only on a single prime. 

For $p$ a prime, $b\in\N_0$ and $n\in\Z$, set $I_{p^b}(n):=1$ if 
$n^2\equiv 4\,(p^{2b})$ and, in case $p=2$, if $(n^2-4)2^{-2b}$ is a discriminant (for $p>2$ this is 
automatically fulfilled). Set $I_{p^b}(n):=0$ otherwise. Define
\begin{equation} \label{12}
  \beta_{(p)}(n) := \sum_{b\ge 0} \frac{1}{p^b} \Big( 1-\frac{1}{p}\,\,\chi_{(n^2-4)p^{-2b}}(p)\Big)^{-1}
  I_{p^b}(n) ,\quad n>2 .
\end{equation}
\begin{lemma} \label{110}
  For $P\ge 2$, we have $\beta_P = \prod_{p\le P} \beta_{(p)}$.
\end{lemma}
\begin{proof}
In the definition of $\beta_P(n)$, write $v=\prod_{p\le P} p^{b_p}$ with $b_p\in\N_0$. There is a 
discriminant $d>0$ with $dv^2=n^2-4$ iff $p^{2b_p}|n^2-4$ for all $p\le P$ and $d:=(n^2-4)v^{-2}\equiv
0,1\,(4)$. Since $p^{2b_p}\equiv 1\,(4)$ for $2<p\le P$, the last condition is equivalent to 
$(n^2-4)2^{-2b_2}\equiv 0,1\,(4)$. If these conditions are fulfilled, we have $(n^2-4)p^{-2b_p} = dr_p^2$
with $r_p\in\N$, $p\not| r_p$, for $p\le P$. Thus $\chi_{(n^2-4)p^{-2b_p}}(p)=\chi_d(p)$ for $p\le P$.
This proves the lemma.
\end{proof}
\begin{cor} \label{109}
  For $q\ge 1$, we have $\alpha\in\calD^q$. In particular, $\lim_{P\to\infty} \beta_P=\alpha$ with respect
  to the $q$-norm.
\end{cor}
\begin{proof}
For $q\in\N$ with $q>1$ fixed and $q'>1$ with $1/(2q)+1/q'=1$ it follows from Lemma~\ref{102} and 
Proposition~\ref{108} that 
\[
  \|\alpha-\beta_P\|_{2q} \ll c_2(P)^{1/q'} + c_3(P)^{1/(2q)} ;
\]
here 
\[
  c_2(P) := \sum_{v>P} \frac{1}{v^{q'}} \to 0 
\]
as $P\to\infty$ since $q'>1$. Furthermore,
\[
  c_3(P) := \sum_{l>P^{2q}} \frac{\tau_{2q}(l)}{l\,K(l)} \to 0
\]
as $P\to\infty$, since
\[
  \sum_{l\ge 1} \frac{\tau_{2q}(l)}{l\,K(l)} = \sum_{a,b\ge 1:\, a \text{ squarefree}} 
  \frac{\tau_{2q}(ab^2)}{ab^2\cdot a} \ll_\eps \sum_{a\ge 1} \frac{a^\eps}{a^2} \sum_{b\ge 1} 
  \frac{b^{2\eps}}{b^2} < \infty .
\] 
Thus $\lim_{P\to\infty} \|\alpha-\beta_P\|_{2q}=0$. For $f:\N\to\infty$ arbitrary and 
$1\le q_1\le q_2<\infty$, we have $\|f\|_{q_1} \le \|f\|_{q_2}$ by H\"older's inequality. Thus 
$\lim_{P\to\infty} \|\alpha-\beta_P\|_q=0$ for all $q\ge 1$. 

Since the $b$-th summand of $\beta_{(p)}$ is $p^{2b+1}$-periodic ($2^{2b+3}$-periodic in case $p=2$) and 
the series representing $\beta_{(p)}$ is uniformly convergent, the function $\beta_{(p)}$ is uniformly
limit periodic, i.e. $\beta_{(p)}\in\calD^u$; here $\calD^u$ is the set of all functions which can be
approximated to an arbitrary accuracy by periodic functions with respect to the supremum norm. Since
$\calD^u$ is closed under multiplication it follows from Lemma~\ref{110} that $\beta_P\in\calD^u$ for 
all $P\ge 2$. This gives $\alpha\in\calD^q$ for all $q\ge 1$.  
\end{proof}
Next the Fourier coefficients of $\alpha$ are computed in terms of the Fourier coefficients of the
$\beta_{(p)}$. In particular, this will show their multiplicativity. 
\begin{lemma} \label{111} \begin{itemize}
    \item[(a)] For all primes $p$, we have $\wh{\beta_{(p)}}(0)=1$.
    \item[(b)] For $b\in\N$, $a\in\Z$, $\gcd(a,b)=1$, choose $a_p\in\Z$ for all $p|b$ such that 
          $\sum_{p|b} a_p p^{-\ord_p b} \equiv ab^{-1} \,(1)$. Then 
          \[
            \wh{\alpha}\big(\frac{a}{b}\big) = \prod_{p|b} \wh{\beta_{(p)}}\Big(
            \frac{a_p}{p^{\ord_p b}}\Big) .
          \]
  \end{itemize}
\end{lemma}
\begin{proof}
From Corollary~\ref{109} it follows that for arbitrary $0<\eps<1$ there is some $P\ge 2$ with 
$\|\alpha-\beta_P\|_1\le\eps$ and $\ord_p b=0$ for all $p>P$. For all $p\le P$ it follows from (\ref{12}) 
that there is some $l_p\ge \ord_p b$ and coefficients $c_p(a_p^\ast)\in\C$, $1\le a_p^\ast\le p^{l_p}$, such
that 
\[
  \Big\| \beta_{(p)} - \sum_{1\le a_p^\ast\le p^{l_p}} c_p(a_p^\ast)\,e_{a_p^\ast/p^{l_p}} \Big\|_u \le 
  \eps' ,
\]
where $\|\cdot\|_u$ denotes the supremum norm and $\eps':=\eps P^{-1} \Big( \max_{p\le P} \|\beta_{(p)}
\|_u+1\Big)^{-P}$. From Lemma~\ref{110} it follows that 
\[
  \Big\| \beta_P - \prod_{p\le P} \Big( \sum_{1\le a_p^\ast\le p^{l_p}} c_p(a_p^\ast)\,e_{a_p^\ast/
  p^{l_p}} \Big) \Big\|_u \le \eps .
\]
Thus 
\[
  \Big\| \alpha -  \sum_{1\le a_p^\ast\le p^{l_p}\,(p\le P)}\,\,  \prod_{p\le P}  c_p(a_p^\ast)\,
  e_{\sum_{p\le P} a_p^\ast/p^{l_p}} \Big\|_1 \le 2\eps .
\]
For $f\in \calD^1$ we have $|\wh{f}|\le \|f\|_1$. Furthermore, $\sum_{p\le P} a_p^\ast p^{-l_p} \equiv
ab^{-1}\,(1)$ iff $a_p^\ast\equiv a_pp^{l_p-\ord_p b}\,(p^{l_p})$ for $p\le P$. Therefore the orthogonality
relation for the exponential function gives 
\[
  \Big| \wh{\alpha}\Big(\frac{a}{b}\Big) - \prod_{p\le P} c_p(a_p p^{l_p-\ord_p b}) \Big| \le 2\eps .
\]
Similarly, $\big|\wh{\beta_{(p)}}(a_pp^{-\ord_p b}) - c_p(a_p p^{l_p-\ord_p b})\big| \le \eps'$ for 
$p\le P$ and thus
\[
  \Big| \prod_{p\le P} \wh{\beta_{(p)}}(a_p p^{-\ord_p b}) - \prod_{p\le P} c_p(a_p p^{l_p-\ord_p b}) \Big|
  \le \eps .
\]
This gives 
\[
  \Big| \wh{\alpha}\Big(\frac{a}{b}\Big) - \prod_{p\le P} \wh{\beta_{(p)}}\Big( \frac{a_p}{p^{\ord_p b}}\Big)
  \Big| \le 3\eps .
\]
In the next section we will compute $\wh{\beta_{(p)}}$ and thereby show that $\wh{\beta_{(p)}}(0)=1$.
This gives (a) and 
\[
  \Big| \wh{\alpha}\Big(\frac{a}{b}\Big) - \prod_{p|b} \wh{\beta_{(p)}}\Big( \frac{a_p}{p^{\ord_p b}}\Big)
  \Big| \le 3\eps .
\]
Since $0<\eps<1$ is arbitrary, (b) follows.
\end{proof}
From Corollary~\ref{109} it follows that (\ref{2}) holds. Here the series on the right hand side is 
absolutely convergent (plug in $r=0$). Thus Lemma~\ref{111} gives 
\[
  \gamma(r) + 1 = \prod_p \Big( 1+\sum_{b_p\ge 1} A_r(p^{b_p}) \Big) ,
\]
where for $p$ prime and $b\in\N$, we define
\[
  A_r(p^b) := \sum_{1\le a\le p^b,\, p\not\,| a} \Big| \wh{\beta_{(p)}}\Big(\frac{a}{p^b}\Big)\Big|^2
  e^{2\pi iar/p^b} .
\]

\section{Computation of the local factors}
\clear

The last step is to calculate $\wh{\beta_{(p)}}$. In particular, this will show that $\wh{\beta_{(p)}}(0)=1$
which is left over from the proof of Lemma~\ref{111}. For $p$ prime, $b\in\N_0$, define
\[
  \beta_{(p,b)}(n) := \Big( 1-\frac{1}{p}\,\,\chi_{(n^2-4)p^{-2b}}(p)\Big)^{-1} I_{p^b}(n) , \quad n>2 .
\]  
Then $\beta_{(p)} = \sum_{b\ge 0} p^{-b}\beta_{(p,b)}$, where the series is uniformly convergent. The
calculation will only be done for $p>2$. The case $p=2$ is similar but somewhat more elaborate. Let 
$b,c\in\N_0$, $a\in\Z$, $p\not| a$. \\
{\bf Case 1:} $2b<c-1$. Since $\beta_{(p,b)}$ is $p^{2b+1}$-periodic, we have $\wh{\beta_{(p,b)}}(a/p^c)=0$.
\\
{\bf Case 2:} $b=c=0$. Then 
\begin{align*}
  &\wh{\beta_{(p,0)}}(0) 
  = \frac{1}{p} \sum_{n\,\bmod\, p} \Big( 1-\frac{1}{p}\,\,\chi_{n^2-4}(p)\Big)^{-1} \\
  &= \frac{1}{p(1-1/p)} \#\big\{ n\bmod p \,\big|\, \chi_{n^2-4}(p)=1 \big\}
     + \frac{1}{p(1+1/p)} \#\big\{ n\bmod p \,\big|\, \chi_{n^2-4}(p)=-1 \big\}  + \frac{2}{p} .
\end{align*} 
The cardinality of the first set is $(p-3)/2$ and that of the second is $(p-1)/2$ (see for example the 
proof of Lemma~3.3 in \cite{peter4}). Thus $\wh{\beta_{(p,0)}}(0)=1-2/(p(p^2-1))$.\\
{\bf Case 3:} $b=0$, $c=1$. Define
\[
  S_p^\pm (a) := \sum_{n\,\bmod\,p:\, \chi_{n^2-4}(p)=\pm1} e^{2\pi ian/p} .
\]
Then 
\begin{align*}
  \wh{\beta_{(p,0)}}\Big(\frac{a}{p}\Big) 
  &= \frac{1}{p} \sum_{n\,\bmod\,p} \Big( 1-\frac{1}{p}\,\,\chi_{n^2-4}(p)
  \Big)^{-1} e^{-2\pi ian/p} \\
  &= \frac{1}{p} \bigg( 2\cos\Big(\frac{4\pi a}{p}\Big) + \sum_\pm \Big(1\mp\frac{1}{p}\Big)^{-1} S_p^\pm(a)
     \bigg) .
\end{align*}
{\bf Case 4:} $2b\ge c-1$, $b\ge 1$. Then 
\[
  \wh{\beta_{(p,b)}}\Big(\frac{a}{p^c}\Big) = \frac{1}{p^{2b+1}} \sum_\pm \sum_{n\,\bmod\,p^{2b+1}:\,
  n\equiv\pm2\,(p^{2b})} \Big( 1-\frac{1}{p}\,\,\chi_{(n^2-4)p^{-2b}}(p)\Big)^{-1} e^{-2\pi ian/p^c} .
\]
Setting $n=\pm 2+mp^{2b}$ gives 
\[
  \wh{\beta_{(p,b)}}\Big(\frac{a}{p^c}\Big) = \frac{1}{p^{2b+1}} \sum_\pm e^{\mp4\pi ia/p^c}
  \sum_{m\,\bmod\,p} \Big( 1-\frac{1}{p}\Big(\frac{m}{p}\Big)\Big)^{-1} e^{\mp 2\pi iamp^{2b-c}} ,
\]
where $\big(\frac{\cdot}{p}\big)$ denotes the Legendre symbol. \\
{\bf Case 4.1:} $2b=c-1$. We have 
\begin{align*}
  \sum_{m\,\bmod\,p:\, (m/p)=1} e^{\mp2\pi i am/p} &= \frac{1}{2} \sum_{m\,\bmod\,p}
  e^{\mp2\pi iam/p} \Big(\Big(\frac{m}{p}\Big)+1\Big) -\frac{1}{2} \\
  &= \frac{1}{2}\sum_{m\,\bmod\,p} e^{\mp 2\pi iam/p}\Big(\frac{m}{p}\Big) - \frac{1}{2} .
\end{align*}
Set $\eps_p:=1$ if $p\equiv1\,(4)$ and $\eps_p:=i$ otherwise. The last sum can be reduced to the Gaussian
sum associated to the Legendre character which can be computed explicitely (see for example 
\cite{davenport}, Chapter~2). This gives for the above quantity the value
\[
  \frac{1}{2} \Big(\frac{\mp a}{p}\Big) \eps_p p^{1/2} - \frac{1}{2} .
\]
Therefore 
\[
  \wh{\beta_{(p,b)}}\Big(\frac{a}{p^c}\Big) = \frac{1}{p^{2b+1}} \bigg( \frac{-2}{p^2-1} \cos\Big(
  \frac{2\pi a}{p^c}\Big) + \frac{p^{3/2}\eps_p}{p^2-1} \Big(\Big(\frac{a}{p}\Big) e^{4\pi i a/p^c}
  +\Big(\frac{-a}{p}\Big) e^{-4\pi ia/p^c}\Big)\bigg) .
\]
{\bf Case 4.2:} $2b>c-1$. Then 
\begin{align*}
  \wh{\beta_{(p,b)}}\Big(\frac{a}{p^c}\Big) &= \frac{1}{p^{2b+1}} \sum_\pm e^{\mp4\pi ia/p^c}
  \sum_{m\,\bmod\,p} \Big( 1-\frac{1}{p}\Big(\frac{m}{p}\Big)\Big)^{-1} \\
  &=  \frac{1}{p^{2b+1}} \sum_\pm e^{\mp4\pi ia/p^c} \Big( \frac{1}{1-1/p}\cdot\frac{p-1}{2} + 
      \frac{1}{1+1/p} \cdot \frac{p-1}{2} + 1 \Big) \\
  &= \frac{2}{p^{2b+1}} \cos\Big(\frac{4\pi a}{p^c}\Big) \frac{p^2+p+1}{p+1} .
\end{align*}

\vspace{0.5cm}
Now we can calculate the Fourier coefficients
\[
  \wh{\beta_{(p)}}\Big(\frac{a}{p^c}\Big) = \sum_{b\ge 0} \frac{1}{p^b} \wh{\beta_{(p,b)}}\Big(
  \frac{a}{p^c}\Big) .
\]
{\bf Case 1:} $c=0$. Then 
\[
  \wh{\beta_{(p)}}(0) = 1-\frac{2}{p(p^2-1)} + \sum_{b\ge 1} \frac{1}{p^b}\frac{2}{p^{2b+1}} \cos(0)
  \frac{p^2+p+1}{p+1} = 1 .
\]
{\bf Case 2:} $c=1$. Since $S_p^+(a)+S_p^-(a)= -2\cos(4\pi a/p)$, we get
\begin{align*}
  \wh{\beta_{(p)}}\Big(\frac{a}{p}\Big) = & \frac{1}{p} \bigg( 2\cos\Big(\frac{4\pi a}{p}\Big)
  +\sum_\pm \Big(1\mp\frac{1}{p} \Big)^{-1} S_p^\pm(a) \bigg) 
  + \sum_{b\ge 1} \frac{1}{p^b} \frac{2}{p^{2b+1}} \cos\Big( \frac{4\pi a}{p}\Big) \frac{p^2+p+1}{p+1} \\
  = & \frac{2p}{p^2-1} \cos\Big(\frac{4\pi a}{p}\Big) + \sum_\pm \frac{1}{p\mp 1} S_p^\pm (a) 
  =  \frac{1}{p^2-1} \sum_{n\,\bmod\,p} \Big(\frac{n^2-4}{p}\Big) e^{2\pi iam/p} .
\end{align*}
{\bf Case 3:} $c\ge 2$. Then 
\begin{align*}
  \wh{\beta_{(p)}}\Big(\frac{a}{p^c}\Big) = & \sum_{0\le b<(c-1)/2} \frac{1}{p^b}\cdot 0 \\
  & + \sum_{b:\, 2b=c-1} \frac{1}{p^b} \frac{1}{p^{2b+1}} \bigg( \frac{-2}{p^2-1} \cos\Big(
    \frac{4\pi a}{p^c}\Big) + \frac{p^{3/2}\eps_p}{p^2-1} \Big(\Big(\frac{a}{p}\Big)e^{4\pi ia/p^c}
    +\Big(\frac{-a}{p}\Big)e^{-4\pi ia/p^c}\Big)\bigg) \\
  & + \sum_{b>(c-1)/2} \frac{1}{p^b} \frac{2}{p^{2b+1}} \cos\Big(\frac{4\pi a}{p^c}\Big) \frac{p^2+p+1}{p+1}
    \\
  = & \frac{1}{(p^2-1)\,p^{3c/2-2}} \left\{\begin{array}{ll} 
      \eps_p\big(\big(\frac{a}{p}\big)e^{4\pi ia/p^c} + \big(\frac{-a}{p}\big) e^{-4\pi ia/p^c}\big) , & 
      c \text{ odd} \\
      2\cos\big(\frac{4\pi a}{p^c} \big) , & c \text{ even} \end{array}\right\} .
\end{align*}

\vspace{0.5cm}
Finally $A_r(p^c)$ can be computed. \\
{\bf Case 1:} $c=1$. Then 
\[
  A_r(p) = \sum_{n_1,n_2\,\bmod\, p} \frac{1}{(p^2-1)^2} \Big(\frac{n_1^2-4}{p}\Big)
  \Big(\frac{n_2^2-4}{p}\Big) \sum_{1\le a\le p-1} e^{2\pi i(n_1-n_2+r)/p} .
\]
The innermost sum is $p-1$ for $n_1-n_2+r\equiv 0\,(p)$ and $-1$ otherwise. Since
\[
  \sum_{n\,\bmod\,p} \Big(\frac{n^2-4}{p}\Big) = \frac{p-3}{2} - \frac{p-1}{2} = -1 ,
\]
the value for $A_r(p)$ follows as given in the theorem.\\
{\bf Case 2:} $c\ge 2$. Then for $p\not| a$, we have
\[
  \Big|\wh{\beta_{(p)}}\Big(\frac{a}{p^c}\Big)\Big| = \frac{1}{(p^2-1)p^{3c/2-2}} 
  \Big| 1+\Big(\frac{-1}{p}\Big)^c\,e^{-8\pi ia/p^c}\Big| .
\] 
Thus 
\[
  A_r(p^c) = \frac{1}{(p^2-1)^2\, p^{3c-4}} \sum_{1\le a\le p^c,\,p\not\,| a} e^{2\pi iar/p^c}
  \Big( 2+2\Big(\frac{-1}{p}\Big)^c \cos\Big(\frac{8\pi a}{p^c}\Big)\Big) . 
\]
A short calculation gives the value in the theorem. \qed

\end{document}